\documentclass[a4paper,12pt]{amsart}
\usepackage{latexsym}
\usepackage{amsmath}
\usepackage{amsfonts}
\usepackage{amssymb}
\usepackage{graphicx}
\newcommand{\R}{\mathbb R}

\newtheorem{theorem}{Theorem}
\newtheorem{lemma}[theorem]{Lemma}
\newtheorem{corollary}[theorem]{Corollary}
\theoremstyle{definition}

\begin{document}
\title[Finitely degenerate hyperbolic problem]{A dyadic decomposition approach to a finitely degenerate hyperbolic problem}
\date{}
\author[M. Cicognani, D. Del Santo]{Massimo Cicognani, Daniele Del Santo}
\address{Massimo  Cicognani, Dipartimento di Matematica,
Universit\`a di Bo\-lo\-gna, Piazza di Porta S. Donato 5, 40127, Bologna, Italy} 
\email{cicognan@dm.unibo.it}
\address{Daniele Del Santo, Dipartimento di Matematica e Informatica,
Universit\`a di Trieste, Via A.~Valerio 12/1, 34127 Trieste, Italy}
\email{delsanto@units.it}
\author[M. Reissig]{Michael Reissig}
\address{Michael Reissig, Institut f\"{u}r Angewandte Analysis, TU
Bergaka\-de\-mie Freiberg, Pr\"uferstrasse 9, 09596 Freiberg, Germany}
\email{reissig@math.tu-freiberg.de}

\dedicatory{To the memory of Stefano Benvenuti}

\begin{abstract}
We use the Littlewood-Paley decomposition technique to obtain a $C^\infty$-well-posedness result for a weakly hyperbolic equation with a finite order of degeneration
\end{abstract}

\maketitle

\def\dif{\partial}
\def\ga{\gamma}
\def\Ga{\Gamma}
\def\ep{\epsilon}
\def\ve{\varepsilon}
\def\La{\Lambda}
\def\la{\lambda}
\def\h{{\bf h}}
\def\be{\beta}
\def\al{\alpha}
\def\de{\delta}
\def\si{\sigma}
\def\lr#1{\langle{#1}\rangle}
\def\qed{$\quad\Box$\vspace{3mm}}
\section{Introduction}
\renewcommand{\theequation}{\arabic{equation}}

Let us consider the following operator
\begin{equation}
L=\partial^2_t -\sum _{j,k=1}^n \partial_{x_j}(a_{jk}(t,x)\partial_{x_k})
+\sum_{j=1}^n b_j(t,x)\partial_{x_j}+c(t,x),
\label{1}
\end{equation}
where $a_{jk}=a_{kj}\in B^\infty([0,T]\times \R^n)$, $b_j$, $c\in C^0([0,T], B^\infty(\R^n))$ 
and $B^\infty$ is the space of the infinitely differentiable functions 
which are bounded with bounded derivatives.
We assume that $L$ is weakly hyperbolic, i.e.
\begin{equation}
a(t,x,\xi)=\sum _{j,k=1}^n a_{jk}(t,x)\xi_j\xi_k/|\xi|^2\geq 0
\label{2}
\end{equation}
for all $(t,x,\xi)\in [0,T]\times \R^n\times (\R^n\setminus\{0\})$. We suppose also that the operator (\ref{1}) has a finite order degeneration in the points in which  $a(t,x,\xi)=0$, i.e. there exists  a positive integer $k$ such that 
\begin{equation}
\sum _{j=0}^k |\partial_t^ja(t,x,\xi)|\not =0
\label{3}
\end{equation}
for all $(t,x,\xi)\in [0,T]\times \R^n\times (\R^n\setminus\{0\})$.
We are interested in the $C^\infty$-well-posedness in the Cauchy problem  for the operator (\ref{1}) with respect to the hyperplane  $\{t=0\}$. It is well-known that a usual hypothesis to such kind of results concerns a relation between the terms of first and second order of the operator, the so-called Levi condition. Our Levi condition reads as follows: that there exist
$C_0>0$, $\gamma\geq 0$ such that 
\begin{equation}
|b(t,x,\xi)|=\left |
\sum _{j=1}^n b_j(t,x)\xi_j/|\xi|\right |\leq C_0 a(t,x,\xi)^\gamma
\label{4}
\end{equation}
for all $(t,x,\xi)\in [0,T]\times \R^n\times (\R^n\setminus\{0\})$.

The $C^\infty$-well-posedness for $L$ satisfying the conditions (\ref{2}), (\ref{3}) and (\ref{4}) with
\begin{equation}
\gamma+{1\over k}\geq {1\over 2} 
\label{5}
\end{equation}
has been proved in \cite{CIO} under the supplementary hypothesis that all the coefficients of the operator $L$ depend only on $t$. A well-known example in \cite{I} shows that the condition (\ref{5}) is sharp. Let us remark that if $k=2$ then condition (\ref{5}) reduces to $\gamma\geq0$, i.e. no Levi condition is requested: it is the case of the effective hyperbolicity studied in \cite{N2}; on the other hand, if $k$ approaches to infinity then
condition (\ref{5}) approaches to $\gamma\geq 1/2$: some special results in this case were obtained in \cite{N}. The result of \cite{CIO} has been extended in \cite{CDFN} to the case that the coefficient $c$ depends also on the variable $x$. The case of  operators $L$ having only time-dependent coefficients in the principal part, while all the other ones depend on $t$ and $x$ has been considered in \cite{CN}, but in this case the 
$C^\infty$-well-posedness has been obtained replacing the condition (\ref{5}) by the more restrictive one
$$
\gamma+{1\over2( k-1)}\geq {1\over 2}  .
$$
Recently, in \cite{AC}, it has been considered the case in which the principal part of the operator $L$ has the following structure:
$$
\partial^2_t-\alpha(t)\sum_{j,k=1}^n\partial_{x_j}(\beta_{jk}(x)\partial_{x_j}).
$$
In this situation the hypotheses (\ref{2}), (\ref{3}), (\ref{4}) and (\ref{5}) are sufficient to prove the $C^\infty$-well-posedness. We remark that the result in \cite{AC} improves that one in \cite{CN}
only in the case of one space variable. 

The proofs of all the above quoted results are based on the use of the so-called approximate energies, which origin goes back to the papers of Colombini, De Giorgi, Jannelli and Spagnolo \cite{CDGS}, \cite{CJS} (a survey on the matter can be found in  \cite{CDS}). In particular in \cite{CIO} the technique is very similar to that one in \cite{CJS}. In \cite{CDFN} and \cite{CN} some convolution weights are introduced to treat the coefficients depending on $x$, while in \cite{AC} the approximate energies method is used together with some pseudodifferential theory arguments.

The approximate energies have been coupled with the Littlewood-Paley decomposition technique in \cite{CL}, where some energy inequalities and some well-posedness results for strictly hyperbolic operators with non-Lipschitz-continuous coefficients have been obtained.  This approach has been successfully used also in the case of strictly hyperbolic operators with oscillating coefficients in \cite{DSKR}. 

In this note we use the approximate energies together with the Littlewood-Paley decomposition technique to improve slightly some of the recalled results of $C^\infty$-well-posedness for weakly hyperbolic  operators with a finite order of degeneration. In particular we will  consider operators having the principal part of the form
$$
\partial^2_t-\alpha(t)\sum_{j,k=1}^n\partial_{x_j}(\beta_{jk}(t,x)\partial_{x_k}),
$$
where 
$$
\sum_{j,k=1}^n\beta_{jk}(t,x)\xi_j\xi_k/|\xi|^2\geq \lambda_0>0
$$
for all $(t,x,\xi)\in [0,T]\times \R^n\times (\R^n\setminus\{0\})$. The proof of this result is inspired by that one of \cite{CIO}, where the use of Fourier transform is replaced by the localization in the phase space given by the dyadic decomposition of Littlewood-Paley.

\bigskip
\noindent
{\bf Acknowledgment}. This work was prepared during the stay of the authors at the Institute of Mathematics of the University of Tsukuba, Japan, in March 2005. The authors would like to express their gratitude to this Institution and in particular to Prof.'s K. Kajitani, T. Kinoshita and S. Wakabayashi.

\section {Main results}

The main result of the paper is contained in the following theorem.

\begin{theorem}
\label{t1}
Let $L$ be the operator defined in (\ref{1}). Suppose 
that the hypotheses (\ref{2}), (\ref{3}),  (\ref{4}) and (\ref{5}) hold.
Suppose moreover that there exist $\alpha$, $\beta_{jk}\in B^\infty$ such that $a_{jk}(t,x)=\alpha(t)\beta_{jk}(t,x)$
for all $j$, $k=1, \dots, n$ 
and for all $(t,x)\in [0,T]\times \R^n$ 
and there exist $\Lambda_0$, $\lambda_0>0$ such that 
\begin{equation}
\Lambda_0\geq \sum_{j,k=1}^n\beta_{jk}(t,x)\xi_j\xi_k/|\xi|^2\geq \lambda_0>0
\label{9}
\end{equation}
for all $(t,x,\xi)\in [0,T]\times \R^n\times (\R^n\setminus\{0\})$.

Then there exist $\delta>0$ such that
\begin{equation}
\begin{array}{l}
\displaystyle{\sup_{0\le t \le T} \left\{\|u(t, \cdot)\|_{H^{m+1-\delta}}  +
 \|\partial_t u(t,\cdot)\|_{H^{m-\delta}}\right\} }\\[0.3 cm]
\quad\displaystyle{ \leq C_m (\|u(0, \cdot)\|_{H^{m+1}}+
 \|\partial_t u(0,\cdot)\|_{H^{m}} + \int_0^T  \| L u(t,\cdot)\|_{
H^{m}}\, dt)}
\end{array}
\label{10}
\end{equation}
for all $m\in \R$ and for all $u\in C^2([0,T], H^\infty(\R^n))$.
\end{theorem}

The usual consequence of the energy inequality (\ref{10}) is stated in the following corollary.

\begin{corollary}
Let $L$ be the operator defined in (\ref{1}). Suppose 
that the hypotheses of the Theorem \ref{t1} hold.

Then the Cauchy problem for $L$ is $C^\infty$-well-posed.
\label{c1}
\end{corollary}

\section{Proofs}

We divide the proof of the Theorem \ref{t1}  into several steps.

\medskip
\noindent
{\it a) the dyadic decomposition}

\medskip
\noindent
We collect in this step some of the well-known facts on the Littlewood-Paley decomposition, referring to \cite{B} and \cite{CL} for the details. Let $\varphi_0\in C^\infty_0(\R^n)$, $0\leq \varphi_0(\xi)\leq 1$, $\varphi_0(\xi)=1$ if $|\xi|\leq 1$, $\varphi_0(\xi)=0$ if $|\xi|\geq 2$, $\varphi_0$ radial and decreasing in $|\xi|$. We set $\varphi(\xi)=\varphi_0(\xi)-\varphi_0(2\xi)$ and, if $\nu$ is an integer greater or equal than 1, $\varphi_\nu(\xi)=\varphi(2^{-\nu}\xi)$. Let $w$ be a function in $H^\infty(\R^n)$; we define
$$
\begin{array}{l}
\displaystyle{w_\nu(x)=\varphi_\nu(D_x)w(x)= {1\over (2\pi)^{n/2}}\int e^{ix\xi}\varphi_\nu(\xi)\hat w(\xi)\, d\xi}\\[0.3 cm]
\qquad\qquad\qquad\qquad\qquad\qquad\quad\displaystyle{={1\over (2\pi)^{n/2}}\int\hat\varphi(2^\nu y)2^{n\nu}w(x-y)\, dy}.
\end{array}
$$
We have that for all $m\in \R$ there exists $K_m>0$ such that
$$
{1\over K_m}\sum_{\nu=1}^{\infty}\|w_\nu\|^2_{L^2} 2^{2m\nu}\leq \|w\|^2_{H^m}\leq K_m \sum_{\nu=1}^{\infty}\|w_\nu\|^2_{L^2} 2^{2m\nu}.
$$
Moreover, for all $\nu\geq 1$, we obtain
\begin{equation}
2^{\nu-1}\|w_\nu\|_{L^2}\leq \|\nabla_x w_\nu\|_{L^2}=(\sum_j \|\partial_{x_j} w_\nu\|^2_{L^2})^{1\over 2}\leq 
 2^{\nu+1}\|w_\nu\|_{L^2}.
\label{13}
\end{equation}

\medskip
\noindent
{\it b) the estimate for the microlocalized approximate energy}

\medskip
\noindent
Let $\varepsilon$ be a positive real number less or equal than 1. Let $u(t,x)$ be a function in $C^1([0,T],H^\infty( \R^n))$.  We set $u_\nu(t,x)=\varphi_\nu(D)u(t,x)$. We introduce the microlocalized approximate energy
$$
\begin{array}{l}
\displaystyle{ 
E_{\nu, \ve}u(t)}\\[0.3 cm]
\quad\displaystyle{ =\|\partial_t u_\nu(t,\cdot)\|^2_{L^2} + \sum^{n}_{j,k=1}\langle(a_{jk}(t,\cdot)+\delta_{jk}\ve)\partial_{x_k}
u_\nu(t,\cdot), \partial_{x_j}u_\nu(t,\cdot) \rangle_{L^2},}
\end{array}
$$
where $\delta_{jk}$ is Kroneker's symbol. An easy computation gives
$$
\begin{array}{l}
\displaystyle{ {d\over dt} E_{\nu, \ve}u(t)}\\[0.3 cm]
\quad\displaystyle{ =2 {\rm Re} \sum_j \langle \ve \partial_{x_j}u_\nu, \partial_{x_j}\partial_t u_\nu\rangle
+ \sum_{j,k} \langle \partial_t a_{jk} \partial_{x_k}u_\nu, \partial_{x_j} u_\nu\rangle}\\[0.3 cm]
\qquad\displaystyle{- 2 {\rm Re} \sum_j \langle b_j\partial_{x_j}u_\nu,\partial_t u_\nu\rangle
- 2 {\rm Re} \langle c u_\nu, \partial_t u_\nu\rangle }\\[0.3 cm]
\qquad\;\;\displaystyle{ +2 {\rm Re} \sum_{j,k} \langle  \partial_{x_j}([\varphi_\nu, a_{jk}]\partial_{x_k}u), \partial_t u_\nu\rangle
+ 2 {\rm Re} \sum_j \langle [\varphi_\nu, b_j] \partial_{x_j}u, \partial_t u_\nu\rangle}\\[0.3 cm]
\qquad\quad\displaystyle{ +2 {\rm Re}  \langle [\varphi_\nu, c]u, \partial_t u_\nu\rangle
+ 2 {\rm Re}  \langle (Lu)_\nu, \partial_t u_\nu\rangle.}\end{array}
$$
We have
$$
\begin{array}{l}
\displaystyle{2 | \sum_j \langle \ve \partial_{x_j}u_\nu, \partial_{x_j}\partial_t u_\nu\rangle|}\\[0.3 cm]
\qquad\qquad\displaystyle{\leq {\ve 2^\nu\over (\alpha(t)+\ve)^{1\over 2} }((\alpha(t)+\ve)(\sum_j\|\partial_{x_j} u_\nu\|^2_{L^2})
+n\|\partial_t u_\nu\|^2_{L^2})}\\[0.3 cm]
\qquad\qquad\displaystyle{\leq C_1 {\ve 2^\nu\over (\alpha(t)+\ve)^{1\over 2}}E_{\nu,\ve}u(t),}
\end{array}
$$
where $C_1>0$ depends only on $\lambda_0$ and $n$; similarly
$$
 \sum_{j,k}| \langle \partial_t a_{jk} \partial_{x_k}u_\nu, \partial_{x_j} u_\nu\rangle|
\leq C_2{|\alpha'(t)|\over \alpha(t)+\ve}E_{\nu,\ve}u(t),
$$
where $C_2>0$ depends only on $\Lambda_0$, $\lambda_0$ and $\sup |\partial_t \beta_{jk}|$. From (\ref{4}) we deduce
that
$$
2 | \sum_j \langle b_j\partial_{x_j}u_\nu,\partial_t u_\nu\rangle|\leq C_3 {1\over (\alpha(t)+\ve)^{{1\over 2}-\gamma}}E_{\nu,\ve}u(t),
$$
where $C_3>0$ depends only on $C_0$ and $\Lambda_0$, and finally,
$$
2 |\langle c u_\nu, \partial_t u_\nu\rangle| \leq C_4 {1\over \ve^{1\over2} 2^\nu}E_{\nu,\ve}u(t),
$$
where again $C_4>0$ depends only on $\sup |c|$. We choose now $\ve=\ve_\nu=2^{-\nu{2k\over 2+k}}$. We remark that with this choice $\ve_\nu^{1\over 2} 2^\nu=2^{\nu{2\over 2+k}}\geq 1$. We obtain 
\begin{equation}
\begin{array}{l}
\displaystyle{ \!\! \!\! \!\!{d\over dt} E_{\nu, \ve_\nu}u(t)}\\[0.3 cm]
\displaystyle{ \leq \tilde C \Big({\ve_\nu 2^\nu\over (\alpha(t)+\ve_\nu)^{1\over 2}}+{|\alpha'(t)|\over \alpha(t)+\ve_\nu}
+{1\over (\alpha(t)+\ve_\nu)^{{1\over 2}-\gamma}}+1\Big)E_{\nu,\ve_\nu}u(t)
 }\\[0.5 cm]
\;\displaystyle{ +2 {\rm Re} \sum_{j,k} \langle  \partial_{x_j}([\varphi_\nu, a_{jk}]\partial_{x_k}u), \partial_t u_\nu\rangle
+ 2 {\rm Re} \sum_j \langle [\varphi_\nu, b_j] \partial_{x_j}u, \partial_t u_\nu\rangle}\\[0.3 cm]
\;\;\displaystyle{ +2 {\rm Re}  \langle [\varphi_\nu, c]u, \partial_t u_\nu\rangle
+ 2 {\rm Re}  \langle (Lu)_\nu, \partial_t u_\nu\rangle},
\end{array}
\label{15}
\end{equation}
where $\tilde C$ depends only on $a_{jk}$, $b_j$ and $c$.

\medskip
\noindent
{\it c) the estimate for the total energy}

\medskip
\noindent
We remind that from \cite[Lemma 1 and 2]{CIO} we have that there exists  $C>0$ depending only on the function $\alpha$ such that
\begin{equation}
\begin{array}{l}
\displaystyle{h(\nu,t)}\\[0.3 cm]
\quad\displaystyle{=\tilde C\int_0^t ({\ve_\nu 2^\nu\over (\alpha(s)+\ve_\nu)^{1\over 2}}+{|\alpha'(s)|\over \alpha(s)+\ve_\nu}
+{1\over (\alpha(s)+\ve_\nu)^{{1\over 2}-\gamma}}+1)\, ds}\\[0.3 cm]
\quad\displaystyle{\leq C\nu}
\end{array}
\label{16}
\end{equation}
for all $t\in [0,T]$ and for all $\nu\geq 1$. We will need 
also the following result; the proof  is very similar 
to those ones for the quoted lemmas from \cite{CIO}: we let it to the reader.
\begin{lemma}
There exists $C>0$ depending only on the function $\alpha$ such that
\begin{equation}
|h(\nu,t)-h(\nu+1, t)|\leq C
\label{17}
\end{equation}
for all $t\in [0,T]$ and for all $\nu\geq 1$.
\label{l1}
\end{lemma}
Now we define the total energy
\begin{equation}
{\tilde E}(t)= \sum_{\nu= 0}^\infty e^{-h(\nu,t)-2\sigma t} E_{\nu, \ve_\nu}u(t),
\label{18}
\end{equation}
where $\sigma$ is a positive constant to be fixed. Our goal is to prove that 
\begin{equation}
{\tilde  E}'(t)\leq \sum_{\nu= 0}^\infty e^{-h(\nu,t)-2\sigma t} \| (Lu)_\nu\|^2_{L^2},
\label{19}
\end{equation}
and from this the inequality  (\ref{10}) will follow. From (\ref{15}) we obtain that
\begin{equation}
\begin{array}{ll}
\displaystyle{{\tilde  E}'(t)\leq}&\displaystyle{
\sum_{\nu= 0}^\infty -2\sigma e^{-h(\nu,t)-2\sigma t} E_{\nu, \ve_\nu}u(t)}\\[0.3cm]
&\displaystyle{+\sum_{\nu= 0}^\infty  e^{-h(\nu,t)-2\sigma t} 2
|\sum_{j,k}\langle [\varphi_\nu, a_{jk}]\partial_{x_k} u, \partial_{x_j}\partial_t u_\nu\rangle|}\\[0.3cm]
&\quad\displaystyle{+\sum_{\nu= 0}^\infty  e^{-h(\nu,t)-2\sigma t} 2
|\sum_{j}\langle [\varphi_\nu, b_j]\partial_{x_j} u, \partial_t u_\nu\rangle|}\\[0.3cm]
&\qquad\displaystyle{+\sum_{\nu= 0}^\infty  e^{-h(\nu,t)-2\sigma t} 2
|\langle [\varphi_\nu, c] u, \partial_t u_\nu\rangle|}\\[0.3cm]
&\qquad\quad\displaystyle{+\sum_{\nu= 0}^\infty  e^{-h(\nu,t)-2\sigma t} 2
|\langle (Lu)_\nu, \partial_t u_\nu\rangle|}.
\end{array}
\label{20}
\end{equation}
We start to estimate 
$$
A=\sum_{\nu= 0}^\infty  e^{-h(\nu,t)-2\sigma t} 2
|\sum_{j,k}\langle [\varphi_\nu, a_{jk}]\partial_{x_k} u, \partial_{x_j}\partial_t u_\nu\rangle|.
$$
We use the technique of \cite[Lemma 4.4]{CL}. We set $\psi_\mu=\varphi_{\mu-1}+
\varphi_\mu+\varphi_{\mu+1}$ ($\varphi_{-1}\equiv 0$); using also (\ref{13})  we have
$$
\begin{array}{ll}
\displaystyle{A}&\leq \displaystyle{\sum_{\nu, \mu= 0}^\infty  e^{-h(\nu,t)-2\sigma t} 2
|\sum_{j,k}\langle([\varphi_\nu, a_{jk}]\psi_\mu)\partial_{x_k} u_\mu, \partial_{x_j}\partial_t u_\nu\rangle|}\\[0.3 cm]
&\displaystyle{\leq \sum_{\nu, \mu= 0}^\infty  e^{-{(h(\nu,t)-h(\mu,t))\over 2}} n2^{\nu+2}
\| [\varphi_\nu, a]\psi_\mu\|_{\mathcal L}}\\[0.3 cm]
&\qquad\qquad\qquad\displaystyle{\cdot e^{-{h(\mu,t)\over 2}-\sigma t}\|\nabla u_\mu\|_{L^2} e^{-{h(\nu,t)\over 2}-\sigma t}\|\partial_t u_\nu\|_{L^2},}
\end{array}
$$
where $\|\cdot\|_{\mathcal L}$ is the norm of bounded linear operators from $(L^2(\R^n_x))^n$ into itself. We remark now that
$$
\begin{array}{ll}
\displaystyle{\| [\varphi_\nu, a]\psi_\mu\|_{\mathcal L}}&\displaystyle{\leq \alpha(t) \| [\varphi_\nu, \beta ]\psi_\mu\|_{\mathcal L}}\\[0.3 cm]
&\displaystyle{\leq(\alpha(t)+\ve_\mu) \| [\varphi_\nu, \beta ]\psi_\mu\|_{\mathcal L}}\\[0.3 cm]
&\displaystyle{\leq c (\alpha(t)+\ve_\mu)^{1\over 2} \| [\varphi_\nu, \beta ]\psi_\mu\|_{\mathcal L},}
\end{array}
$$
where $c$  depends only  on $\alpha$. Finally
$$
(\alpha(t)+\ve_\mu)^{1\over 2}\|\nabla u_\mu\|_{L^2}
\leq c' |\sum_{j,k}\langle (a_{jk} +\delta_{jk}\ve_\mu)\partial_{x_k} u_\mu, 
\partial_{x_j} u_\mu\rangle|^{1\over 2},
$$
where $c'$ depends only on $\beta_{jk}$. We deduce that 
$$
\begin{array}{ll}
\displaystyle{A}&\leq \displaystyle{c''  \sum_{\nu, \mu= 0}^\infty  e^{-{(h(\nu,t)-h(\mu,t))\over 2}} 2^{\nu}
\| [\varphi_\nu, a]\psi_\mu\|_{\mathcal L}}\\[0.5 cm]
&\quad\displaystyle{\cdot e^{-{h(\mu,t)\over 2}-\sigma t}|\sum_{j,k}\langle (a_{jk} +\delta_{jk}\ve_\mu)\partial_{x_k} u_\mu, 
\partial_{x_j} u_\mu\rangle|^{1\over 2}
 e^{-{h(\nu,t)\over 2}-\sigma t}\|\partial_t u_\nu\|_{L^2}.}
\end{array}
$$
The idea is now to use Schur's lemma. We introduce
$$
k_{\nu \mu}(t) =e^{-{(h(\nu,t)-h(\mu,t))\over 2}} 2^{\nu} \| [\varphi_\nu, a]\psi_\mu\|_{\mathcal L}.
$$
We have to consider
$$
\sup_{\nu} \sum_{\mu} |k_{\nu \mu}(t)|, \qquad
\sup_{\mu} \sum_{\nu} |k_{\nu \mu}(t)|.
$$
To this end the following lemma will be useful.
\begin{lemma}
Let $|\nu-\mu|\leq 2$.

Then there exists  $C$ depending only on $\beta$ such that
\begin{equation}
\| [\varphi_\nu, \beta]\psi_\mu\|_{\mathcal L}\leq C 2^{-\nu}.
\label{21}
\end{equation}
 Let $|\nu-\mu|\geq 3$. 
 
 Then for all $N>0$ there exists $C_N>0$ depending only on $N$ and $\beta$ such that
\begin{equation}
\| [\varphi_\nu, \beta]\psi_\mu\|_{\mathcal L}\leq C_N 2^{-N\cdot \max \{\nu, \mu\}}.
\label{22}
\end{equation}
\label{l2}
\end{lemma}

\smallskip
\noindent
{\it Proof.}
It is easy to prove (\ref{21}) (see \cite[Prop. 3.6]{CL}). The inequality (\ref{22}) can be proved taking into account the fact that the supports of $\varphi_\nu$ and $\psi_\mu$ are disjoint and using the asymptotic formula of  of $ [\varphi_\nu, \beta]\psi_\mu$
(see \cite[Prop. 4.5]{CL}).\hfill\qed

We denote now by $\tilde k_{\nu \mu}(t)$ the value $ k_{\nu \mu}(t)\chi_{\{|\nu-\mu|\leq 2\}}(\nu,\mu)$, where $\chi_\Omega$ is the characteristic function of the set $\Omega$ . Then
$$
\sum_{\mu} |\tilde k_{\nu \mu}(t)|=\sum_{j=\nu-2}^{\nu+2}  | k_{\nu j}(t)|.
$$
From (\ref{17})  and  (\ref{21}) we have that $\sum_{\mu} |\tilde k_{\nu \mu}(t)|\leq C$ for all $t\in [0,T]$, where $C$ does not depend on $\nu$. Similarly $\sum_{\nu} |\tilde k_{\nu \mu}(t)|\leq C$ for all $t\in [0,T]$ and for all $\mu$. 
We denote by $ k^\ast_{\nu \mu}(t)$ the value $ k_{\nu \mu}(t)\chi_{\{|\nu-\mu|\geq 3\}}(\nu,\mu)$.
From (\ref{16}) we have that 
$$
|h(\nu,t)-h(\mu,t)|\leq C_1(\nu+\mu)
$$
for all $\nu$, $\mu$ and for all $t\in [0,T]$. Consequently using (\ref{22}) we deduce
$$
| k^\ast_{\nu \mu}(t)|\leq C_N  e^{C_1(\nu+\mu) +\nu-N\cdot \max\{\nu,\mu\}}.
$$
Easily we have
$$
\sum_\mu | k^\ast_{\nu \mu}(t)|\leq C,
$$
where $C$ does not depend on $\nu$; similarly in the remaining case. We obtain
$$
\begin{array}{ll}
\displaystyle{A}&\leq \displaystyle{C \Big(\sum_{\nu= 0}^\infty   e^{-{h(\mu,t)\over 2}-\sigma t}|\sum_{j,k}\langle (a_{jk} +\delta_{jk}\ve_\nu)\partial_{x_j} u_\nu, 
\partial_{x_k} u_\nu\rangle|\Big)^{1\over 2}}\\[0.5 cm]
 & \qquad\qquad\qquad\qquad\qquad\qquad\qquad\displaystyle{\cdot  \Big(\sum_{\nu= 0}^\infty 
 e^{-{h(\nu,t)\over 2}-\sigma t}\|\partial_t u_\nu\|^2_{L^2}\Big)^{1\over 2},}
\end{array}
$$
so that there exists $C>0$ which does not depend on $\sigma$ such that
$$
A\leq C\tilde E(t).
$$
Let us consider briefly the next term in (\ref{20}). We set 
$$
B=\sum_{\nu= 0}^\infty   e^{-{h(\nu,t)}-2\sigma t}2|\sum_j\langle [\varphi_\nu, b_j]\partial_{x_j}u, \partial_t u_\nu\rangle|.
$$
As before 
$$
\begin{array}{ll}
\displaystyle{B}
&\leq \displaystyle{ \sum_{\nu, \mu= 0}^\infty  e^{-{(h(\nu,t)-h(\mu,t))\over 2}} 2
\| [\varphi_\nu, b]\psi_\mu\|_{\mathcal L}}\\[0.5 cm]
&\qquad\qquad\qquad \displaystyle{\cdot
e^{-{h(\mu,t)\over 2}-\sigma t}\|\nabla u_\mu\|_{L^2} e^{-{h(\nu,t)\over 2}-\sigma t}\|\partial_t u_\nu\|_{L^2}}\\[0.5 cm]
&\leq\displaystyle{c\sum_{\nu, \mu= 0}^\infty  e^{-{(h(\nu,t)-h(\mu,t))\over 2}} 
\| [\varphi_\nu, b]\psi_\mu\|_{\mathcal L}\ve_\mu^{-1}}\\[0.5 cm]
&\quad \displaystyle{\cdot e^{-{h(\mu,t)\over 2}-\sigma t}|\sum_{j,k}\langle (a_{jk} +\delta_{jk}\ve_\mu)\partial_{x_j} u_\mu, 
\partial_{x_k} u_\mu\rangle|\;
e^{-{h(\nu,t)\over 2}-\sigma t}\|\partial_t u_\nu\|_{L^2}.}
\end{array}
$$
A computation similar to the previous one gives
$$
B\leq C\tilde E(t).
$$
The estimate of the other terms from (\ref{20}) is strightforward. We finally obtain
$$
\tilde E'(t)\leq (-2\sigma +C) \tilde E(t) +\sum_{\nu=0}^\infty  e^{-{h(\nu,t)}-2\sigma t}\|(Lu)_\nu\|^2_{L^2},
$$
where $C$ depends on $a_{jk}$, $b_j$, $c$ but does not depend on $\sigma$. We choose $\sigma>C/2$ and (\ref{19}) follows. From this the inequality (\ref{10}) is obtained with a standard argument (see \cite[p. 695]{CL}).  The proof is complete.

\end{document}